\input amstex
\magnification=\magstep0
\documentstyle{amsppt}
\pagewidth{6.0in}
\input amstex
\topmatter
\title
A note on a question due to A. Garsia
\endtitle
\author
Tewodros Amdeberhan\\
May 03, 2009\\
Tamdeber\@tulane.edu
\endauthor
\abstract
Adriano Garsia has provided an explicit formula to enumerate a certain class of permutations in the symmetric group $S_n$. In this short note, we prove a stronger version of the assertion under a specific conjecture. 
\endabstract
\endtopmatter
\def\({\left(}
\def\){\right)}

\document
\noindent
\smallskip
\noindent
\bf 1. Introduction and nomenclature \rm
\bigskip
\noindent
In this section, we present a question due to A. Garsia [1] and we lay down the relevant framework as well as some basic results. In the sequel, assume $\binom{c}d=0$ whenever $c<d$ or $d<0$. We start with the actual problem proposed by Garsia. It is stated as follows.
\bigskip
\noindent
\bf Problem G. \it Fix $k\in\Bbb{Z}_{+}$. For each $n\geq 2k$, define the set of permutations
$$A_{n,n-k}:=\{\mu=a_1a_2\cdots a_n\in S_n: a_1<a_2<\cdots<a_{n-k}; \text{no incr. subseq. of length $>n-k$}\}.$$
Prove (by elementary means) the enumeration
$$\#A_{n,n-k}:=
\sum_{i=0}^k(-1)^{k-i}\binom{k}i\frac{n!}{(n-i)!}.\tag1$$ \rm
Before discussing our approach let us introduce a few notations. Consider the set of permutations in $A_{n,n-k}$ that begin with (prefix) $i\in[1,n]$, and is denoted by, 
$$B_{n,n-k}(i):=\{\mu\in A_{n,n-k}: a_1=i\}.$$
Notice that if $i>k+1$ then $B_{n,n-k}(i)$ is empty. On such account these sets are disregarded and we will restrict $1\leq i\leq k+1$. Observe that $A_{n,n-k}$ is now a disjoint union of the $B_{n,n-k}(i)$ and thus
$$\#A_{n,n-k}=\sum_{i=1}^{k+1}\#B_{n,n-k}(i).\tag2$$

\noindent
Let $T$ stand for transposing a matrix. Call the column vector $$\tilde{\bold{B}}_{n,n-k}:=[\#B_{n,n-k}(1),\#B_{n,n-k}(2),\dots,\#B_{n,n-k}(k+1)]^T,$$
associated with $\#A_{n,n-k}$, to be the \it component vector. \rm 

\smallskip
\noindent
For $k$ fixed, we shall always begin with $n=2k$. This motivates us to form the vector, which we call the \it $k$-kernel, \rm
$$\bold{K}:=\tilde{\bold{B}}_{2k,k}=[\#B_{2k,k}(1),\#B_{2k,k}(2),\dots,\#B_{2k,k}(k+1)]^{T}.$$ 
We are ready to state and prove the following recursive relation.
\smallskip
\noindent
\bf Lemma 1.1 \it For each $i\in[1,k+1]$ and $n\geq 2k$, it holds 
$$\#B_{n+1,n+1-k}(i)=\sum_{r=i}^{k+1}\#B_{n,n-k}(r).\tag3$$ \rm
\bf Proof. \rm Suppose $\mu=a_1a_2\cdots a_n\in B_{n,n-k}(r)$. Then the map $\sigma_{r,i}: B_{n,n-k}(r)\rightarrow B_{n+1,n+1-k}(i)$ defined by $\sigma_{r,i}: \mu\mapsto i\cdot\nu_i(a_1)\nu_i(a_2)\cdots\nu_i(a_n)$ where
$$\nu_i(a_j)=\cases 1+a_j \qquad\text{if $a_j\geq i$}\\
a_j \qquad{} \qquad\text{otherwise},\endcases$$
is injective. And $\sigma_{r,i}$ is a bijection when the elements of $B_{n+1,n+1-k}(i)$ are partitioned by their second entries. $\square$
\smallskip
\noindent
\bf Corollary 1.2 \it For each $i\in[1,k+1]$ and $n\geq 2k$, it holds that
$\tilde{\bold{B}}_{n,n-k}=\bold{C}_{n,n-k}\bold{K}$ where
$$\bold{C}_{n,n-k}:=\left(\binom{r+n-2k-i-1}{r-i}\right)_{i,r=1}^{k+1};$$
or, in detail, 
$$\#B_{n,n-k}(i)=\sum_{r=i}^{k+1}\binom{r+n-2k-i-1}{r-i}\#B_{2k,k}(r).$$ \rm
\bf Proof. \rm This is immediate from a repeated application of Lemma 1.1. $\square$
\bigskip
\noindent
\bf 2. Illustrative Examples \rm
\bigskip
\noindent
The next three tabular examples elaborate how equation (3) enables the successive generation of new columns and thereby determining $\#A_{n,n-k}$ via equation (2).
\smallskip
\noindent
\bf Example 2.1. \rm Let $k=1, n\geq 2$. Then
$$\pmatrix
\#B_{n,n-1}(i)&n=2&n=3&n=4&n=5&n=6\\
\hdotsfor 6\\
\#B_{n,n-1}(1)&0&1&2&3&4\\
\#B_{n,n-1}(2)&1&1&1&1&1\\
\hdotsfor 6\\
\#A_{n,n-1}=&\bold1&\bold2&\bold3&\bold4&\bold5
\endpmatrix$$
\bf Example 2.2. \rm Let $k=2, n\geq 4$. Then
$$\pmatrix
\#B_{n,n-2}(i)&n=4&n=5&n=6&n=7&n=8&n=9\\
\hdotsfor 7\\
\#B_{n,n-2}(1)&1&5&11&19&29&41\\
\#B_{n,n-2}(2)&2&4&6&8&10&12\\
\#B_{n,n-2}(3)&2&2&2&2&2&2\\
\hdotsfor 7\\
\#A_{n,n-2}=&\bold5&\bold{11}&\bold{19}&\bold{29}&\bold{41}&\bold{55}
\endpmatrix$$
\bf Example 2.3. \rm Let $k=3, n\geq 6$. Then
$$\pmatrix
\#B_{n,n-3}(i)&n=6&n=7&n=8&n=9&n=10&n=11\\
\hdotsfor 7\\
\#B_{n,n-3}(1)&14&47&104&191&314&479\\
\#B_{n,n-3}(2)&15&33&57&87&123&165\\
\#B_{n,n-3}(3)&12&18&24&30&36&42\\
\#B_{n,n-3}(4)&6&6&6&6&6&6\\
\hdotsfor 7\\
\#A_{n,n-3}=&\bold{47}&\bold{104}&\bold{191}&\bold{314}&\bold{479}&\bold{692}
\endpmatrix$$
\bf Example 2.4. \rm We demonstrate the proof of Lemma 1.1 when $n=4$ and $k=2$. The tables show a listing of the permutations in the sets $B_{\alpha,\beta}(\gamma)$. Beginning with and based on
$$A_{4,4-2}=
\pmatrix
B_{4,4-2}(1)&B_{4,4-2}(2)&B_{4,4-2}(3)\\
\hdotsfor 3\\
\bold{14}32&\bold{24}13&\bold{34}12\\
{}&\bold{24}31&\bold{34}21
\endpmatrix$$
we construct (dashed lines indicate empty contributions, bold numbers are monotonic)
$$\align
B_{4+1,4+1-2}(1)=
\pmatrix
\bold{125}43&\bold{135}24&\bold{145}23\\
-&\bold{135}42&\bold{145}32
\endpmatrix,
\qquad
&B_{4+1,4+1-2}(2)=
\pmatrix
-&\bold{235}14&\bold{245}13\\
-&\bold{235}41&\bold{245}31
\endpmatrix,\\
&B_{4+1,4+1-2}(3)=
\pmatrix
-&-&\bold{345}12\\
-&-&\bold{345}21
\endpmatrix.\endalign$$
\bf 3. Main Results and a Conjecture \rm
\bigskip
\noindent
The conclusion of Corollary 1.2 says that if the $k$-kernel $\bold{K}$ is known then $\#A_{n,n-k}$ as well as its component vector $\tilde{\bold{B}}_{n,n-k}$ can be computed, for any $n\geq 2k$. Thus, the main task is \it how to determine the kernel vector $\bold{K}$. \rm Although at present we do not have a proof, we are convinced that the conjecture given below addresses the question fully. 
\bigskip
\noindent
Let us define a column vector $\bold{V}_k:=[v_1,v_2,\dots,v_{k+1}]^T$ that we call the \it initial $k$-vector \rm according to
$$v_i:=\sum_{b\geq 0}(-1)^{k-b}\binom{k}b\binom{i}bb!,$$
and the $(k+1)\times(k+1)$-matrix $\bold{M}_k$ by
$$\bold{M}_k:=\left(\binom{i+j-2k-1}{i-2k}\right)_{i,j=1}^{k+1}.$$ \rm
\bf Conjecture 3.1 \it The $k$-kernel vector $\bold{K}=[\#B_{2k,k}(j): 1\leq j\leq k+1]^T$ equals to $\bold{K}=\bold{M}_k^{-1}\bold{V}_k$. \rm
\bigskip
\noindent
The next result is new. It offers a stronger statement than Problem G in allowing us to calculate any component vector $\tilde{\bold{B}}_{n,n-k}$, for all $n\geq 2k$. 
\bigskip
\noindent
\bf Lemma 3.2. \it Suppose Conjecture 3.1 is true and let $n\geq 2k$. Then $\tilde{\bold{B}}_{n,n-k}=\bold{Q}_{k,n}^{-1}\bold{V}_k$ where
$$\bold{Q}_{k,n}:=\left((-1)^{j-1}\binom{n-i-1}{j-1}\right)_{i,j=1}^{k+1}.$$ \rm
\bf Proof. \rm From Cor. 1.2 and Conj. 3.1, the claim amounts to $\bold{Q}\bold{C}=\bold{M}$. See Appendix. $\square$
\bigskip
\noindent
We now supply a restatement of Lemma 3.2 that, in some sense, avoids inverting a matrix. 
\bigskip
\noindent
\bf Proposition 3.3. \it Both matrices $\bold{M}_k$ and $\bold{Q}_{k,n}$ have determinant $1$, hence $\bold{M}_k^{-1}$ and $\bold{Q}^{-1}_{n,k}$ consist of integral entries. In particular, the $i$-th entry of $\tilde{\bold{B}}_{n,n-k}$ equals to $\#B_{n,n-k}(i)=\det(\bold{Q}_{n,k}\vert_i \bold{V}_k)$, i.e. the determinant of the matrix $\bold{Q}_{n,k}$ whose $i$-th column is replaced by $\bold{V}_k$. Similarly the $i$-th entry of $\bold{K}$ equals $\#B_{2k,k}(i)=\det(\bold{M}_k\vert_i\bold{V}_k)$. \rm
\bigskip
\noindent
\bf Proof. \rm We generalize the matrices by adding free parameters. For instance, set $x=-2k$ and $y=-1$ to recover $\bold{M}_k$.
$$\det\left(\binom{i+j+x+y}{i+x}\right)_{i,j=1}^{k+1}=
\prod_{i=1}^{k+1}\binom{i+1+x+y}{i+x}\binom{i+y}{1+y}^{-1}.$$
This new determinant is perfectly amenable to \it Dodgson condensation. \rm $\square$
\bigskip
\noindent
\bf Proof of Problem G. \rm Let $\bold{1}_k$ stand for the $(k+1)$-row vector $[1,1,\dots,1]$. From Lemma 3.2 and the definition of our vectors, we obtain  
$$\#A_{n,n-k}=\bold{1}_k\cdot\tilde{\bold{B}}_{n,n-k}=
\bold{1}_k\cdot\bold{Q}_{k,n}^{-1}\bold{V}_k.$$
Next, we solve the equation $\bold{U}_{k,n}\bold{Q}_{k,n}=\bold{1}_k$ for a row vector $\bold{U}_{k,n}$. It turn out that $\bold{U}_{k,n}=[u_1,\dots,u_{k+1}]$ where
$$u_j=(-1)^{k+1-j}\binom{n-2}k\binom{k}{j-1}\frac{n-1}{n-j}.$$
The verification involves routine binomial identities (see Appendix). The final step too revolves around identites; namely,
$$\align
\#A_{n,n-k}&=\bold{U}_{k,n}\bold{V}_k\\
&=\sum_{j=1}^{k+1}u_j\sum_{b\geq0}(-1)^{k-b}\binom{k}b\binom{j}bb!\\
&=(n-1)\binom{n-2}k\sum_{b\geq0}(-1)^b\binom{k}bb!\sum_{j=1}^{k+1}
\frac{(-1)^{j-1}}{n-j}\binom{k}{j-1}\binom{j}b\\
&=(n-1)\binom{n-2}k\sum_{b\geq0}(-1)^b\binom{k}bb!\left\{\frac{(-1)^k}{n-1}\frac{\binom{n}b}{\binom{n-2}k}\right\}\\
&=\sum_{b\geq0}(-1)^{k-b}\binom{k}b\binom{n}bb!.
\endalign$$
We have thus arrived at the desired formula, hence the proof is complete. $\square$
\bigskip
\noindent
\bf Problem. \it Find a combinatorial proof of Lemma 3.2 for the components vector. \rm

\pagebreak

\noindent
\bf APPENDIX \rm
\bigskip
\noindent
In this section, we append some identities that were needed in proving the results of the previous sections. There are several ways to achieve this but we just mention that the justifications can be carried out using the automatic method of Wilf and Zeilberger. 
\bigskip
\noindent
\bf Lemma A. \it We have $\bold{Q}\bold{C}=\bold{K}$. \rm
\bigskip
\noindent
\bf Proof. \rm Since $(-1)^c\binom{c}d=\binom{-c-1+d}d$, we convert $(-1)^{r-1}\binom{n-i-1}{j-1}=\binom{r-n+i-1}{r-1}$. Combining this with the definition of the corresponding matrices, the assertion tantamount to the identity
$$\sum_{r=1}^j\binom{r-n+i-1}{r-1}\binom{j+n-2k-r-1}{j-r}=\binom{i+j-2k-1}{j-1}.$$
This, however, is the special case $y=r-1, x=j-1, A=i-n, B=n-2k-1$ of the \it Vandermonde-Chu \rm convolution formula
$$\sum_{y=0}^x\binom{y+A}{y}\binom{x-y+B}{x-y}=\binom{A+B+x+1}{x}.\qquad \square$$
\bf Lemma B. \rm We have $\bold{U}_{k,n}\bold{Q}_{k,n}=\bold{1}_k$. \rm
\bigskip
\noindent
\bf Proof. \rm The claim is equivalent to the identity 
$$\Phi(k,r):=\sum_{j=1}^{k+1}(-1)^{k+r+j}\frac{(n-1)}{n-j}\binom{n-2}k\binom{k}{j-1}\binom{n-j-1}{r-1}=1;
\qquad \text{where $1\leq r\leq k+1$}.$$
It is interesting to note that this identity is invalid for $r>k+1$, so we should be a bit careful in our analysis. This is achieved by extracting two recurrences using \it Zeilberger algorithm: \rm
$$\align
(k-r+2)(n-k-2)&\{\Phi(k+1,r)-\Phi(k,r)\}-
(k+2)(n-k-3)\{\Phi(k+2,r)-\Phi(k+1,r)\}=0,\\
r(n-r-1)&\{\Phi(k,r+1)-\Phi(k,r)\}-
(k-r)(r+1)\{\Phi(k,r+2)-\Phi(k,r+1)\}=0.\endalign$$
The proof is completed by global induction on $k$, and an internal induction on $r$. $\square$
\bigskip
\noindent
\bf Lemma C. \it For $0\leq b\leq k$, we have $\sum_{j=1}^{k+1}\frac{(-1)^{j-1}}{n-j}\binom{k}{j-1}\binom{j}b=\frac{(-1)^k}{n-1}\frac{\binom{n}b}{\binom{n-2}k}$. \rm
\bigskip
\noindent
\bf Proof. \rm As a first step, re-write the identity at hand in the form
$$\psi(k,b):=\sum_{j=1}^{k+1}(-1)^{k+j-1}\frac{(n-1)\binom{n-2}k}{(n-j)\binom{n}b}\binom{k}{j-1}\binom{j}b=1.$$ \rm
Once more, this formula works for any $n$ and $k$ but only for $0\leq b\leq k$. Zeilberger algorithm yields two recursive relations. We find first order recurrences for both parameters $k$ and $b$; that is, $\psi(k+1,b)-\psi(k,b)=0$ and $\psi(k,b+1)-\psi(k,b)=0$. Now proceed as in Lemma B. $\square$


\bigskip


\Refs
\widestnumber\key{1} 

\ref \key 1 \by R. P. Stanley \paper \it Private communication 
\endref

\endRefs
\enddocument